\documentclass{amsart}
\usepackage{enumerate}
\usepackage{amssymb,latexsym}
\usepackage{graphicx}
\usepackage{verbatim}
\usepackage{xcolor} 
\usepackage{hyperref}
\hypersetup{
    linktoc=all,     
    linkcolor=black,  
}

\newtheorem{theorem}{Theorem}
\newtheorem{corollary}[theorem]{Corollary}
\newtheorem{lemma}[theorem]{Lemma}
\newtheorem{proposition}[theorem]{Proposition}

\newcommand{\R}{\mathbb{R}}

\newcommand{\ce}{C_{\epsilon}}
\newcommand{\se}{S_{\epsilon}}

\begin{document}

\title{Cyclic conformally flat hypersurfaces revisited }

\maketitle
\begin{center}
\author{Jo\~ao Paulo dos Santos        \and
        Ruy Tojeiro}  
      \footnote{The first author is supported by FAPDF 0193.001346/2016. The second author is partially 
supported by Fapesp grant 2016/23746-6 and 
CNPq grant 303002/2017-4.}
\end{center}
\date{}

\begin{abstract}
In this article we classify the conformally flat Euclidean hypersurfaces of dimension three with three distinct principal curvatures  of $\mathbb{R}^4$, $\mathbb{S}^3\times \mathbb{R}$ and $\mathbb{H}^3\times \mathbb{R}$ with the property that the tangent component of the vector field $\partial/\partial t$ is a principal direction at any point. Here  $\partial/\partial t$ stands  for either a constant unit vector field in $\mathbb{R}^4$ or the unit vector field tangent to the factor $\mathbb{R}$ in the product spaces $\mathbb{S}^3\times \mathbb{R}$ and 
$\mathbb{H}^3\times \mathbb{R}$, respectively. Then we use this result to give a simple proof  of an alternative classification of the  cyclic conformally flat hypersurfaces of $\mathbb{R}^4$, that is, the conformally flat hypersurfaces of $\mathbb{R}^4$ with three distinct principal curvatures such that the curvature lines correspondent to one of its principal curvatures are extrinsic circles. We also characterize the cyclic conformally flat hypersurfaces of $\mathbb{R}^4$ as those conformally flat hypersurfaces of dimension three with three distinct principal curvatures for which there exists a conformal Killing vector field of $\mathbb{R}^4$ whose tangent component is an eigenvector field correspondent to one of its principal curvatures.  \vspace{2ex}

\end{abstract}

\noindent \emph{2020 Mathematics Subject Classification:} 53 B25.\vspace{2ex}

\noindent \emph{Key words and phrases:} {\small {\em conformally flat hypersurface,  linear Weingarten surface, \\
 product spaces, conformal Killing vector field.}}


\section{Introduction}

  E. Cartan proved in \cite{ca} that if $f\colon M^{n} \to \mathbb{R}^{n+1}$  is an isometric immersion of a Riemannian manifold $M^n$ of dimension   $n\geq 4$, then $M^n$ is conformally flat  if and only if  $f$  has a principal curvature of multiplicity 
at least $n-1$. Recall that a Riemannian manifold $M^n$ is conformally flat if each point of $M^n$ has an open neighborhood that is conformally diffeomorphic to an open subset of Euclidean space $\R^n$. Cartan also proved that any hypersurface $f\colon M^{3} \to \mathbb{R}^{4}$ with a principal curvature of multiplicity greater than one is conformally flat, and  realized that the converse is no longer true in this case.
Thus, generic conformally flat Euclidean hypersurfaces of   dimension   $n\geq 4$ are envelopes of one-parameter families of hyperspheres, but in dimension  $n=3$ there appears an interesting further class of conformally flat hypersurfaces which have three distinct principal curvatures.

  Cartan's investigations were taken up by Hertrich-Jeromin \cite{hj}, who showed  that any conformally flat Euclidean hypersurface of dimension three with three distinct principal curvatures carries local principal coordinates $u_1,u_2,u_3$ with respect to which the induced metric can be written as
$$
ds^2=\sum_{i=1}^3 v_i^2du_i^2,
$$
with the Lam\'e coefficients $v_i$, $1\leq i\leq 3$, satisfying the Guichard condition, say, $v_2^2=v_1^2+v_3^2$.
Then he used the conformal invariance of this condition to associate with
each such hypersurface a  Guichard net in $\R^3$, that is, a conformally
flat metric on an open subset of $\R^3$ satisfying the Guichard condition,
which is unique up to a Moebius transformation. He also proved in \cite{hj}
(see also Section 2.4.6 in \cite{hj2}) that each conformally flat $3$-metric 
satisfying the Guichard condition gives rise to a unique 
(up to a Moebius transformation) conformally flat hypersurface in $\R^4$ 
(see also \cite{ct1}).

 Improving earlier work by Suyama (see \cite{su1}, \cite{su2}), 
Hertrich--Jeromin and Suyama \cite{hjs} gave a classification of conformally flat 
 hypersurfaces  whose associated  Guichard nets in $\R^3$ are cyclic, 
 that is, one of their coordinate line families consist of circular arcs. These include the so-called \emph{conformal product} conformally flat hypersurfaces, which are either cylinders over umbilic-free surfaces of constant Gauss curvature in $\mathbb{R}^3$, cones over  umbilic-free surfaces of constant Gauss curvature in $\mathbb{S}^3$ or rotation hypersurfaces over  umbilic-free surfaces in a half-space $\R_+^3$ of $\R^3$, regarded as the half-space model of $\mathbb{H}^3$,  which have  constant Gauss curvature with respect to the metric induced from the hyperbolic metric on $\R_+^3$. Conformal product conformally flat hypersurfaces in $\R^4$ have been characterized in \cite{dt1} as those conformally flat hypersurfaces with three distinct principal curvatures in $\R^4$  such that the curvature lines correspondent to one of its principal curvatures are arcs of circles or straight lines in $\mathbb{R}^4$.
 
  A class of  noncyclic conformally flat hypersurfaces was subsequently studied in \cite{hjs2}, whose associated  Guichard systems are of Bianchi-type, that is,  its coordinate surfaces have constant Gauss curvature. However, until not very long ago, all the known explicit examples of conformally flat hypersurfaces of $\R^4$ with three distinct principal curvatures belonged to the class of
cyclic conformally flat hypersurfaces. More recently, a Ribaucour transformation for the class of conformally flat hypersurfaces of $\R^4$ with three distinct principal curvatures, based on the characterization of such hypersurfaces obtained in \cite{ct1}, was developed in \cite{ct2}, which allowed to construct explicit noncyclic examples  (see also \cite{hjsuy} and \cite{st}).\vspace{1ex}    

 One of the goals of this article is to give a simple proof of an 
 alternative description of cyclic conformally flat hypersurfaces of $\mathbb{R}^4$. 
 This will be derived as a consequence of a classification of independent interest
 of the conformally flat hypersurfaces of dimension three with three distinct principal curvatures  of $\mathbb{R}^4$, $\mathbb{S}^3\times \mathbb{R}$ and $\mathbb{H}^3\times \mathbb{R}$ with the property that the tangent component of $\partial/\partial t$ is a principal direction at any point. Here  $\partial/\partial t$ stands  for either a constant unit vector field in $\mathbb{R}^4$ or the unit vector field tangent to the factor $\mathbb{R}$ in $\mathbb{S}^3\times \mathbb{R}$ and 
$\mathbb{H}^3\times \mathbb{R}$.
 
 First recall that  there exists a conformal diffeomorphism 
 $\Phi\colon  \mathbb{S}^{n-1}\times\mathbb{R} \to \mathbb{R}^n \setminus\{0\}$ 
 given by $(x,t) \mapsto e^t x$.
Similarly, there is  a  conformal diffeomorphism 
$$
\Psi \colon \mathbb{H}^n\times \mathbb{S}^{1}\subset \mathbb{R}^{n+1}_1 
\times \mathbb{R}^{2} \to \mathbb{R}^{n+1}\setminus \mathbb{R}^{n-1}
$$
onto the complement of a subspace $\R^{n-1}\subset \R^{n+1}$ given as follows. 
Choose a  pseudo-orthonormal 
basis $e_0, e_1, \ldots, e_{n-1},e_n$  of the Lorentzian space $\R^{n+1}_1$ with
$\left<e_0,e_0\right>=0=\left<e_n, e_n\right>$, 
$\left<e_0, e_n\right>=-1/2$ and $\left<e_i,e_j\right>=\delta_{ij}$
for $1\leq i\leq n-1$ and $0\leq j\leq n$.  Then 
\[
\Psi(x_0 e_0 + \ldots + x_{n} e_n, (y_1, y_2)) = \frac{1}{x_0}(x_1, 
\ldots, x_{n-1}, y_1, y_2).
\]
Composing $\Psi$ with the isometric covering map 
$$
\pi\colon \mathbb{H}^n\times \mathbb{R}\to \mathbb{H}^n\times \mathbb{S}^1: (x,t) \mapsto (x, (\cos t, \sin t))
$$ 
produces a conformal covering map 
$ \Phi \colon \mathbb{H}^n\times \mathbb{R}\to \mathbb{R}^{n+1}\setminus \mathbb{R}^{n-1}$
given by
\begin{equation}\label{tpsi}
\Phi(x_0 e_0 + \ldots + x_n e_n, t) = 
\frac{1}{x_0}(x_1, \ldots, x_{n-1}, \cos t, \sin t).
\end{equation}

   In what follows,  $\mathbb{Q}_\epsilon^3 \subset \R_\mu^{3+|\epsilon|}$ denotes $\mathbb{S}^3$ if $\epsilon=1$, $\mathbb{R}^3$ if $\epsilon=0$  and  $\mathbb{H}^3$ if $\epsilon=-1$, with $\mu=0$ if $\epsilon=0$ or $\epsilon=1$, and $\mu=1$ if $\epsilon=-1$. Given a surface  $h\colon\, M^{2}\to \mathbb{Q}_\epsilon^3$, let $h_s \colon\, M^{2}\to \mathbb{Q}_\epsilon^3\subset \R_\mu^{3+|\epsilon|}$ be the the family of its parallel surfaces, that is,   
$$h_s(x)=C_\epsilon(s)h(x)+S_\epsilon(s)N(x),$$
where $N$ is a unit normal vector field to $h$ and the functions $\ce$ and $\se$ are given by
$$
C_\epsilon(s)=\left\{\begin{array}{l}
\cos s, \,\,\,\mbox{if}\,\,\epsilon=1
\vspace{1.5ex}\\
1, \,\,\,\mbox{if}\,\,\epsilon=0
\vspace{1.5ex}\\
\cosh s, \,\,\,\mbox{if}\,\,\epsilon=-1
\end{array}\right.\,\,\,\,\,\,\,\,\mbox{and}\,\,\,\,\,\,\,\,\,\,\,
S_\epsilon(s)=\left\{\begin{array}{l}
\sin s, \,\,\,\mbox{if}\,\,\epsilon=1
\vspace{1.5ex}\\
s, \,\,\,\mbox{if}\,\,\epsilon=0
\vspace{1.5ex}\\
\sinh s, \,\,\,\mbox{if}\,\,\epsilon=-1.
\end{array}\right.
$$

The classification of the conformally flat hypersurfaces of dimension three with three distinct principal curvatures  of $\mathbb{R}^4$, $\mathbb{S}^3\times \mathbb{R}$ and $\mathbb{H}^3\times \mathbb{R}$ with the property that the tangent component of $\partial/\partial t$ is a principal direction at any point is as follows.

\begin{theorem} \label{thm:main2} Let $h\colon\, M^{2}\to \mathbb{Q}_\epsilon^3$ be an umbilic-free linear Weingarten surface, i.e., the extrinsic curvature $K_{ext}$ and the mean curvature
 $H$ of $h$ satisfy
\begin{equation}\label{linear-weingarten-F}
PK_{ext}+Q H = R 
\end{equation}
for some $P$, $Q$, $R\in \mathbb{R}$.
 Set $\overline{P}=P+\epsilon R$, $\overline{Q}=Q$, $\overline{R}=P-\epsilon R+4$ and $\Lambda = 2(\epsilon^2-1) R $, let  $I\subset \mathbb{R}$  be an open interval  where
 $$
 0<r(s):=\frac{1}{4}(\overline{P}\ce(2s) +\overline{Q} \se(2s) + \Lambda \se^2(s) + \overline{R})<1
 $$ 
 and let $a\colon I\to \mathbb{R}$ be the smooth function on $I$ given by
 \begin{equation} \label{function-a}
a(s)=\int_{s_0}^s\sqrt{ \frac{1-r(s)}{r(s)}}\;ds,\;\;s_0\in I.
\end{equation}
Then the map
$f\colon M^{2}\times I\to \mathbb{Q}_\epsilon^3\times \R\subset \R_\mu^{4+|\epsilon|}$ given by
\begin{equation}\label{eq:constantangle}f(x,s)=h_s(x)+a(s)\partial/\partial t,
\end{equation}
where $\partial/\partial t$ denotes a unit vector field tangent to  $\mathbb{R}$, defines, on the open subset $ M^3\subset M^{2}\times I$ of its regular points, a conformally flat hypersurface with three distinct principal curvatures such that the tangent component of 
 $\partial/\partial t$ is a principal direction of $f$ at any point. 

 Conversely, any conformally flat hypersurface $f\colon M^3\to \mathbb{Q}_\epsilon^3\times \mathbb{R}$ with these properties is given in this way. 
\end{theorem}

  Regarding $h$ as an isometric immersion into $\R^{4+|\epsilon|}_\mu$, its normal space at each point $x\in M^2$ is a vector space whose dimension is $2$ if $\epsilon=0$ and $3$ otherwise, and which is either Lorentzian or Riemannian, according to whether $\mu=1$ or $\mu=0$, respectively. If $\epsilon\neq 0$, it is spanned by the position vector $h(x)$, the normal vector $N(x)$ to $h$ in $\mathbb{Q}^3_\epsilon$ at $x$ and the constant vector $\partial/\partial t$. Notice that these give rise to parallel vector fields along $h$ with respect to its normal connection. For a fixed $x\in M^2$, we can regard $s\mapsto f(x, s) = \ce(s)h(x) + \se(s)N(x) + a(s)\partial/\partial t$, when $\epsilon\neq 0$, as a curve in a
cylinder $\mathbb{Q}_\epsilon^1\times \mathbb{R}$, with axis $\partial/\partial t$, contained in the normal space of $h$ at $x$. Thus $f(M)$ is generated by parallel transporting such curve along $h$ with respect to its normal connection.\vspace{1ex}

 We derive from Theorem \ref{thm:main2} the following alternative classification of cyclic conformally flat hypersurfaces of $\mathbb{R}^4$.
 In the next statement, the map $\Phi$ denotes either the conformal diffeomorphism 
 $\Phi\colon \mathbb{Q}_\epsilon^3\times \R\to \R^{4}\setminus \{0\}$ if $\epsilon=1$, the conformal covering map $\Phi\colon \mathbb{Q}_{\epsilon}^3\times \R\to \R^{4}\setminus \R^2$ if $\epsilon=-1$ or the isometry $\Phi\colon \mathbb{Q}_\epsilon^3\times \R\to  \R^4$ if $\epsilon=0$.

\begin{theorem} \label{thm:main} Let $f\colon M^3 \to \mathbb{Q}_\epsilon^3\times \mathbb{R}$ be a hypersurface given as in Theorem~\ref{thm:main2}. Then $F=\mathcal{I} \circ \Phi \circ f\colon M^3 \to \mathbb{R}^4$, where  $\mathcal{I}$ is either the identity map or an inversion with respect to a hypersphere  in $\R^4$, is a cyclic conformally flat hypersurface.

Conversely, any cyclic conformally flat hypersurface 
of $\mathbb{R}^4$ is given in this way. 
\end{theorem}

Conformal product conformally flat hypersurfaces correspond to the case in which the hypersurface $f\colon M^3 \to \mathbb{Q}_\epsilon^3\times \mathbb{R}$  in Theorem~\ref{thm:main2} is a vertical cylinder over an umbilic-free surface with constant Gauss curvature  $h\colon M^2\to\mathbb{Q}_\epsilon^3$. Alternatively, in terms of the geometric description after Theorem~\ref{thm:main2}, $f(M)$ is generated by parallel transporting, with respect to the normal connection of $h$, a straight line in the direction of $\partial/\partial t$ in a fixed normal space of $h$ in $\R^{4+|\epsilon|}_\mu$.\vspace{1ex}

 In the classification of cyclic conformally flat hypersurfaces given in \cite{hjs}, the authors deal with hypersurfaces in the sphere $\mathbb{S}^4$, taking into account the invariance of the conditions involved under conformal diffeomorphisms between the ambient space forms. Then they use a Moebius geometric technology to show that any such hypersurface can be produced, up to such a conformal diffeomorphism, from a hypersurface in some space form that is given in terms of a family of parallel linear Weingarten surfaces and a solution of a certain pendulum-type ordinary differential equation. Our approach is somewhat more elementary in nature, and the parametrization of cyclic conformally flat hypersurfaces $f\colon M^3\to \mathbb{R}^4$ provided  by Theorems \eqref{thm:main2} and \eqref{thm:main} only requires  a single integration.\vspace{1ex}

 Theorem \ref{thm:main} also yields the following characterization of cyclic conformally flat hypersurfaces. Let $x_1, \ldots, x_{n+1}$ denote  the standard coordinates in $\R^{n+1}$ and let $\partial_{x_i}$ be a unit vector field tangent to the $x_i$-coordinate curve, $1\leq i\leq n+1$. It is well known that the Lie algebra of conformal Killing vector fields in $\R^{n+1}$ has dimension $\frac{1}{2}(n+2)(n+3)$
and is generated by the  constant vector fields $\partial_{x_i}$, $1\leq i\leq n+1$, the Killing vector fields $\mathcal{K}_{ij}=x_i\partial_{x_j}- x_j\partial_{x_i}$, $1\leq i\neq j\leq n+1$, generating rotations around the linear subspaces $\R^{n-1}$ of $\R^{n+1}$
given by  $x_{i} = 0 = x_{j}$, by the radial vector field
$\mathcal{R}=\displaystyle{\sum_{i=1}^{n+1}x_i \partial_{x_i}}$,
and by the vector fields
$$
\mathcal{C}_i=\frac{1}{2}(x_i^2-\sum_{j \neq i}x_j^2)\partial_{x_i}+x_i\sum_{j \neq i}x_j\partial_{x_j},\;\;1\leq i\leq n+1.
$$
 
 \begin{corollary} \label{cor:ckilling} A conformally flat hypersurface  
$f\colon M^3\to \mathbb{R}^4$  with three distinct principal curvatures is cyclic if and only if the  tangent component of one of the above conformal Killing vector fields is an eigenvector field correspondent to one of its principal curvatures.  
\end{corollary}

\section{Proof of  Theorem \ref{thm:main2}}

It is well known that a three-dimensional  Riemannian manifold $M^3$ is conformally flat if and only if its Schouten tensor $L=T-(3/2)sI$ satisfies the Codazzi equation
\begin{equation}\label{eq:codazzi}
(\nabla_X L)Y=(\nabla_Y L)X
\end{equation}
for all $X, Y\in \mathfrak{X}(M)$ (see, e.g, \cite{dt2}, p. 545), where $(\nabla_X L)Y=\nabla_X LY- L(\nabla_X Y)$. Here $T$ is the endomorphism associated with the Ricci tensor 
and $s$ is the scalar curvature.

  Assume that $M^3$ carries local coordinates $x_1, x_2, x_3$ with respect to which its Riemannian metric can be written as
  $$g=v_1^2dx_1^2+v_2^2dx_2^2+v_3^2dx_3^2,$$
   where $v_1$, $v_2$ and $v_3$ are smooth functions. In the sequel, we denote by $\psi_i$ the partial derivative of a function $\psi$ with respect to $x_i$ and  by $\psi_{ij}$ its second order partial derivative $\partial_{x_i} \partial_{x_j} \psi$. Denote
   \begin{equation}\label{eq:phiij} {\phi}^{ij}=\frac{{{v}_{j,i}}}{{v}_i},\;\;1\leq i\neq j\leq 3,
   \end{equation}
where ${{v}_{j,i}}$ denotes the derivative of $v_j$ with respect to $x_i$. Let $\partial_1, \partial_2, \partial_3$ be the coordinate vector fields and set $X_k=v_k^{-1} \partial_k$, $1\leq k\leq 3$.
Then the curvature tensor $R$ of $g$ satisfies (see \cite{dt2}, p. 20)
   \begin{equation}\label{eq:R1}
    R(\partial_i, \partial_j)X_k=\left( ({\phi}^{kj})_i-{\phi}^{ki}{\phi}^{ij}\right)X_j-\left(({\phi}^{ki})_j-{\phi}^{kj}{\phi}^{ji}\right)X_i,\;\;1\leq i\neq j\neq k\neq i\leq 3,
    \end{equation}
   and
  $$
   -v_iv_jK_{ij}=-\left\langle R(\partial_i, \partial_j)X_j, X_i\right\rangle=   
    ({\phi}^{ij})_i+ ({\phi}^{ji})_j
   +\phi^{ki}\phi^{kj}, \;\;1\leq i\neq j\neq k\neq i\leq 3,
   $$  
where $K_{ij}$, $1\leq i\neq j\leq 3$, is the sectional curvature along the plane spanned by $\partial_i$ and $\partial_j$. Now suppose further that 
$$v_1=e^\alpha,\;\;\;v_2=e^\beta\;\;\;\mbox{and}\;\;\;v_3\equiv 1$$
for some smooth functions $\alpha$ and $\beta$ satisfying 
\begin{equation}\label{f23-h13}
\alpha_{23}+\alpha_2(\alpha-\beta)_3 = 0\;\;\;\mbox{and}\;\;\;
\beta_{13}-\beta_1(\alpha-\beta)_3 = 0.
\end{equation}
Notice that the preceding equations are equivalent to 
$ ({\phi}^{ij})_3=0$ for $1\leq i\neq j\leq 2$. Then Eq. \eqref{eq:R1} implies that
$R(\partial_3, \partial_i)X_j=0$ for $1\leq i\neq j\leq 2$, and it follows that
 $L\partial_i=\ell_i\partial_i$, $1\leq i\leq 3$,
where 
$$2\ell_1=K_{12}+K_{13}-K_{23}, \;\;2\ell_2=K_{12}+K_{23}-K_{13}\;\;\mbox{and}\;\;
   2\ell_3=K_{13}+K_{23}-K_{12}.$$
    
Define $\psi_i=\ell_i{v}_i$, $1\leq i\leq 3$. Then $L$ satisfies the Codazzi equation \eqref{eq:codazzi} if and only if 
$$
\psi_{i,j}={\phi}^{ji} \psi_j,\;\;\;1\leq i\neq j\leq 3.
$$        
  Then,  a  straightforward computation yields the following lemma.
  
\begin{lemma}\label{le:suyama} The metric
\begin{equation}\label{eq:tildeg} \tilde{g} = e^{2\alpha} dx^2_1 + e^{2\beta} dx_2^2 + dx_3^2,
\end{equation}
with the smooth functions $\alpha$ and $\beta$ satisfying \eqref{f23-h13}, 
 is conformally flat if and only  $\alpha$ and $\beta$ satisfy the partial differential equations
 \begin{equation}
\begin{array}{c}
(e^{-2\beta}(\alpha_{22}+(\alpha_{2})^{2}-\alpha_{2}\beta_{2}))_{2}+(e^{-2\alpha}(\beta_{11}+(\beta_{1})^{2}-\alpha_{1}\beta_{1}))_{2} \vspace{1ex}\\
-(\alpha_{33}+(\alpha_{3})^{2}+\beta_{33}+(\beta_{3})^{2}-\alpha_{3}\beta_{3})_{2} = 0,
\end{array} \label{conf4-new}
\end{equation}

\begin{equation}
\begin{array}{c}
(e^{-2\beta}(\alpha_{22}+(\alpha_{2})^{2}-\alpha_{2}\beta_{2}))_{1}+(e^{-2\alpha}(\beta_{11}+(\beta_{1})^{2}-\alpha_{1}\beta_{1}))_{1}\vspace{1ex} \\
-(\alpha_{33}+(\alpha_{3})^{2}+\beta_{33}+(\beta_{3})^{2}-\alpha_{3}\beta_{3})_{1} = 0
\end{array} \label{conf5-new}
\end{equation}
and
\begin{equation}
\begin{array}{c}
e^{-2\beta}(\alpha_{22}+(\alpha_{2})^{2}-\alpha_{2}\beta_{2})_{3}+(e^{-2\alpha}(\beta_{11}+(\beta_{1})^{2}-\alpha_1\beta_{1}))_{3} \vspace{1ex}\\
+(\alpha_{3}\beta_{3}+\beta_{33}+(\beta_{3})^{2}-\alpha_{33}-(\alpha_{3})^{2})_{3} \vspace{1ex}\\
=2(\alpha_{33}+(\alpha_{3})^{2}-\alpha_{3}\beta_{3})\beta_{3} -2e^{-2\alpha}(\beta_{11}+(\beta_{1})^{2}-\alpha_1\beta_{1})\beta_{3}.
\end{array} \label{conf6-new}
\end{equation}
\end{lemma}
\vspace*{1ex}

\noindent \emph{Proof of  Theorem \ref{thm:main2}:} By Theorem $1$ in \cite{to2}, if $f\colon M^{2}\times I\to \mathbb{Q}_\epsilon^3\times \mathbb{R}$ is given by (\ref{eq:constantangle}) in terms of an arbitrary  surface $h\colon\, M^{2}\to \mathbb{Q}_\epsilon^3$ and a smooth function $a\colon I\to \mathbb{R}$, then the tangent component of 
 $\partial/\partial t$ is a principal direction at any point of the restriction of $f$ to the subset $M^3\subset M^{2}\times I$  of its regular points. Conversely, any hypersurface  $f\colon M^{3}\to \mathbb{Q}_\epsilon^3\times \mathbb{R}$ with this property is given in this way. Therefore, we must prove that $f$ has three distinct principal curvatures and the  metric induced by $f$ on $M^3$ is conformally flat if and only if $h$ and $a$ are as in the statement.

   The metric induced by $f$ on $M^3$ is
$$
d\sigma^2= b^2(s)ds^2+g_s,
$$
where $b(s)=\sqrt{1+(a'(s))^2}$ and $g_s$ is the metric induced by $h_s$. 

 For some of the computations that follow, it is convenient to have in mind the following relations between the functions $\ce$ and $\se$:
$$
\begin{array}{l}
\ce^2(s) + \epsilon \se^2(s) = 1, \\
\ce(2s) = \ce^2(s) - \epsilon \se^2(s)  \,\, \textnormal{ and } \,\, \se(2s) = 2 \ce(s) \se(s) , \\
\ce'(s) = - \epsilon \se(s) \,\, \textnormal{ and } \,\, \se'(s) = \ce(s).
\end{array}
$$

Let $N$ be a unit normal vector field to $h$ and let $N_s$ be the unit normal vector field to $h_s$ given by $N_s(x)=C_\epsilon N(x)-\epsilon S_\epsilon h(x)$. Then 
$$\eta(x,s)=-\frac{a'(s)}{b(s}N_s(x)+\frac{1}{b(s)}\frac{\partial}{\partial t}$$
defines a unit normal vector field to $f$, and the shape operators $A$ of 
$f$ at  $(x,s)$ and  $A^s$ of $h_s$ at $x$ with respect to $\eta$ and $N_s$, respectively, are related by 
\begin{equation}\label{eq:sff}
A X = -\dfrac{a'(s)}{b(s)} A^s X, \, \textnormal{ if } X \in T_xM^2, \;\;\mbox{and}\;\;
A\partial_s=\frac{a''(s)}{b^3(s)}\partial_s,
\end{equation}
where $\partial_s$ is a unit vector field tangent to $I$. Thus $f$ has three distinct principal curvatures at  $(x,s)$ if and only if $x$ is not an umbilic point for  $h_s$, and hence for $h$. 

Under the assumption  that $h$ has no umbilic points, there exist  locally principal  coordinates $x_1, x_2$ on $M_2$ with respect to which the first and second fundamental forms  of $h$ are  
$$I={v}_1^2dx_1^2+ {v}_2^2dx_2^2\;\;\;\mbox{and}\;\;\;II=V_1{v}_1dx_1^2+ V_2{v}_2dx_2^2,$$
respectively. Therefore the first fundamental form of $h_s$ and its second fundamental form with respect to $N_s$ are given, respectively, by
$$I^s=(v^s_1)^2dx_1^2+ (v^s_2)^2dx_2^2\;\;\;\mbox{and}\;\;\;II^s=V^s_1{v}^s_1dx_1^2+ V^s_2{v}^s_2dx_2^2,$$
where
$${v}^s_i=C_\epsilon v_i - S_\epsilon V_i\;\;\;\mbox{and}\;\;\;
V_i^s=\epsilon S_\epsilon v_i+C_\epsilon V_i=- v_{i,3}^s, \;\; 1\leq i\leq 2.$$
Notice that for all $s\in I$ we have 
\begin{equation}\label{eq:inds}
\frac{ v_{j,i}^s}{v_i^s}=\frac{v_{j,i}}{v_i}, \;\;1\leq i\neq j\leq 2,
\end{equation}
and that ${\displaystyle k_i^s=\frac{V_i^s}{v_i^s}}$, $1\leq i\leq 2$, are the principal curvatures of $h_s$. It follows easily that the extrinsic curvature and mean curvature
\begin{equation}\label{eq:kshs}
K^s_{ext}=k_1^sk_2^s=\frac{V_1^sV_2^s}{v_1^sv_2^s}\;\;\mbox{and}\;\;H^s=\frac{1}{2}\left(\frac{V_1^s}{v_1^s}+\frac{V_2^s}{v_2^s}\right)
\end{equation}
of $h_s$ are related to the extrinsic curvature $K_{ext}$ and the mean curvature $H$ of $h$ by
$$
\begin{array}{rcl}
K^s_{ext}&=&\dfrac{\epsilon^2\se^2(s)+\epsilon\se(2s)H+\ce^2(s)K_{ext}}{\ce^2(s)-\se(2s)H+\se^2(s)K_{ext}}, \vspace{1ex} \\
H^s&=&\dfrac{\epsilon\se(2s)+2\ce(2s)H-\se(2s)K_{ext}}{2\left(\ce^2(s)-\se(2s)H+\se^2(s)K_{ext}\right)}.
\end{array}
$$

  Since $\partial_s$  is a principal direction of $f$ by the second equation in \eqref{eq:sff}, then  $x_1, x_2, x_3:=s$ are local principal coordinates for $f$ with respect to which its induced metric is given by $b^2\tilde g$, where $\tilde g$ has the form \eqref{eq:tildeg} with 
  \begin{equation}\label{eq:efeh}
  e^\alpha=\frac{v_1^s}{b}\;\;\;\;\mbox{and}\;\;\;\; e^\beta=\frac{v_2^s}{b}.
  \end{equation}
  
  It follows from \eqref{eq:inds} that the functions $\phi^{ij}$,  associated with the metric $\tilde g$ by means of \eqref{eq:phiij}, satisfy $(\phi^{ij})_3=0$ for $1\leq i\neq j\leq 2$. Thus $\alpha$ and $\beta$ satisfy \eqref{f23-h13}. We now investigate when they also satisfy (\ref{conf4-new}), (\ref{conf5-new}) and (\ref{conf6-new}).

In terms of the function $\rho$ defined by
\begin{equation}
\begin{array}{rcl}
\rho &:=& e^{-2\beta}(\alpha_{22}+(\alpha_{2})^{2}-\alpha_{2}\beta_{2})+e^{-2\alpha}(\beta_{11}+(\beta_{1})^{2}-\alpha_{1}\beta_{1}),\vspace{1ex} \\
&& -\left(\alpha_{33}+(\alpha_{3})^{2}+\beta_{33}+(\beta_{3})^{2}-\alpha_{3}\beta_{3}\right).
\end{array} \label{rho}
\end{equation}
Eqs. (\ref{conf4-new}), (\ref{conf5-new}) and (\ref{conf6-new}) are equivalent to
\begin{equation}
\begin{array}{rcl}
\rho_i &=& 0, \, \, 1\leq i\leq 2, \vspace{1ex}\\
\rho_3 &=& -2\beta_3 \rho - 2\beta_3 (\beta_{33}+\beta_3^2)-2 (\beta_{33}+(\beta_3)^2)_3
\end{array} \label{partial-derivatives-rho}
\end{equation}

 Differentiating (\ref{eq:efeh}) we obtain 
  $$\alpha_2 = \phi^{21}e^{\beta-\alpha}\;\;\;\mbox{and}\;\;\; \beta_1 = \phi^{12}e^{\alpha-\beta},$$
  hence
  \begin{equation}\label{h11f22}
\beta_{11}+\beta_1(\beta-\alpha)_1 = \phi^{12}_1 e^{\alpha-\beta}\;\;\;\mbox{and}\;\;\; \alpha_{22}+\alpha_2(\alpha-\beta)_2 = \phi^{21}_2 e^{\beta-\alpha}. 
\end{equation}
Therefore, the sum of the first two terms on the right hand-side of \eqref{rho} is
$$
e^{-(\alpha+\beta)}\left(\phi_2^{21}+\phi_1^{12}\right)
= -\frac{b^2}{v_1^sv_2^s}(V_1^sV_2^s+\epsilon v_1^sv_2^s)=-b^2(K^s_{ext}+\epsilon),
$$
where the first equality  follows from the Gauss equation of $h_s$, bearing in mind that $\phi^{ij}=\phi^{ij}_s$ for $1\leq i\neq j\leq 2$. 
  
On the other hand, setting $B=\log b$ we have
\begin{equation}
\begin{array}{rcl}
\alpha_3 &=& - \left( \dfrac{V_1^s}{v_1^s} + \dfrac{b'}{b} \right)= - ( k_1^s + B'),\vspace{1ex} \\
\alpha_{33}&=& - \left( \epsilon + \left( \dfrac{V_1^s}{v_1^s} \right)^2 + \dfrac{b''}{b} - \left(\dfrac{b'}{b}\right)^2  \right)=- ( \epsilon + (k_1^s)^2 + B''), \vspace{1ex} \\
\beta_3 &=& - \left( \dfrac{V_2^s}{v_2^s} + \dfrac{b'}{b} \right)= - ( k_2^s + B'), \label{derivatives-f-h-x3} \vspace{1ex}\\
\beta_{33}&=& - \left( \epsilon + \left( \dfrac{V_2^s}{v_2^s} \right)^2 + \dfrac{b''}{b} - \left(\dfrac{b'}{b}\right)^2  \right)=- ( \epsilon + (k_2^s)^2 + B'').
\end{array}
\end{equation}

 It follows from  \eqref{derivatives-f-h-x3} that
 \begin{equation}\label{f33h33}
\alpha_{33} + (\alpha_3)^2 + \beta_{33} + (\beta_3)^2 - \beta_3 \alpha_3 =  2B'H^s-K^s_{ext}-2\epsilon -2B'' +(B')^2.
\end{equation}
From \eqref{rho}, \eqref{h11f22} and \eqref{f33h33} we obtain
$
\rho = 2 B'' - (B')^2 - \epsilon e^{2B} + 2 \varepsilon - \varphi, 
$
where 
\begin{equation}\label{weingarten-part1}\varphi(x_1, x_2, x_3) := \left( e^{2B} - 1 \right)K_{ext}^s + 2 B' H^s.
\end{equation}
 Using the last two equations in \eqref{derivatives-f-h-x3}, the equations in \eqref{partial-derivatives-rho} are equivalent to
\begin{equation}
\begin{array}{rcl}
\varphi_i &=&0, \, 1 \leq i \leq 2, \vspace{1ex} \\
\varphi_3 &=& 2 k_2^s \theta + 2B'(\epsilon+\varphi), \label{F-prime-1}
\end{array}
\end{equation}
where 
\begin{equation}\label{eq:theta}
\theta= \epsilon e^{2B} + B''-\epsilon + \varphi - 2 (B')^2 .
\end{equation} 

 We have shown so far that, for a hypersurface $f\colon M^{2}\times I\to \mathbb{Q}_\epsilon^3\times \mathbb{R}$ given by (\ref{eq:constantangle}) in terms of 
an umbilic-free  surface $h\colon\, M^{2}\to \mathbb{Q}_\epsilon^3$ and a smooth function $a\colon I\to \mathbb{R}$, the  metric induced by $f$ on the subset $M^3\subset M^{2}\times I$ of its regular points  is conformally flat if and only if the preceding equations are satisfied. We now show that this is the case if and only if $h$ and $a$ are as in the statement.

   Suppose first that Eqs. \eqref{F-prime-1} hold. The first two equations imply that $\varphi$ depends only on $x_3=s$. Let us compute the derivative of $\varphi$ with respect to $x_3$. Differentiating \eqref{eq:kshs} with respect to $x_3$ by using that
$v_{i,3}^s=-V_i^s$ and $V_{i,3}^s=\epsilon v_i$, $1\leq i\leq 2$,  gives 
$$(K^s_{ext})_3 = 2 H^s(\epsilon+K_{ext}^s )\;\;\;\mbox{and}\;\;\;(H^s)_3 = \epsilon + 2 (H^s)^2 - K_{ext}^s.$$
The preceding relations yield
\begin{equation}
\varphi_3 = 2 H^s \theta + 2 B'(\epsilon + \varphi). \label{F-prime-2}
\end{equation}
Comparing \eqref{F-prime-1} and \eqref{F-prime-2} gives
$
0 = (k_1^s - k_2^s) \theta.
$
Since $k_1^s \neq k_2^s$ at any point, we conclude that $\theta$ is identically zero.
Therefore, the functions $\varphi$ and $B$ satisfy the system of ordinary differential equations
\begin{equation}
\left\{
\begin{array}{rcl}
\varphi' - 2B'(\epsilon+\varphi) &=& 0, \\
\epsilon + 2(B')^2 - \epsilon e^{2B} - B'' - \varphi &=& 0.
\end{array} \label{system-F-B}
\right.
\end{equation}
The first equation of \eqref{system-F-B} implies that
\begin{equation}\label{eq:varphi}
\varphi = \lambda e^{2B}-\epsilon
\end{equation}
for some  $\lambda\in \mathbb{R}$. Substituting this formula into the second one, it becomes
$$
(e^{-2B})''+ 4 \epsilon e^{-2B} = 2 (\epsilon + \lambda),
$$
whose general solution is
\begin{equation}\label{eq:B}
e^{-2B(x_3)} = c_1 \ce(2x_3) + c_2 \se(2x_3) + (1-\epsilon^2) \lambda \se^2(x_3) +  \dfrac{\epsilon(\epsilon+\lambda)}{2}
\end{equation}
for some $c_1, c_2\in \mathbb{R}$. Since $B=\log b$, the above equation is equivalent to
\begin{equation}
(1+(a')^2)^{-1} = c_1 \ce(2x_3) + c_2 \se(2x_3) + (1-\epsilon^2) \lambda \se^2(x_3) + \dfrac{\epsilon(\epsilon+\lambda)}{2}. \label{ode-for-a}
\end{equation}

Evaluating \eqref{weingarten-part1} at $x_3=0$ by using \eqref{eq:varphi} and \eqref{eq:B} implies that  
$$
\left[\epsilon(\epsilon+\lambda)+2c_1-2 \right]K_{ext}+(4c_2)H = \epsilon^2(\epsilon+\lambda)+2\epsilon c_1 - 2 \lambda.
\label{linear-weingarten-F2}
$$
Set 
$$
P=\epsilon(\epsilon+\lambda)+2c_1-2,\;\;\;Q=4c_2\;\;\;\mbox{and}\;\;\;R=\epsilon^2(\epsilon+\lambda)+2\epsilon c_1 - 2 \lambda.
$$

 Since $c_1$, $c_2$ and $\lambda\in \mathbb{R}$ are arbitrary, we see that also $P, Q, R\in \mathbb{R}$ are arbitrary, and now \eqref{ode-for-a} implies that $a$ is given as in the statement.

Conversely, let $h\colon\, M^{2}\to \mathbb{Q}_\epsilon^3$ be an umbilic-free linear Weingarten surface, let $a\colon I\to \mathbb{R}$ be a  smooth function on the open interval $I\subset \mathbb{R}$ as in the statement and let $f$ be given by \eqref{eq:constantangle}. We must show that the function $\varphi$, defined by \eqref{weingarten-part1} on the open subset $M^3\subset M^2\times I$ of regular points of $f$, with $B=\log b$ and $b=(1+(a')^2)^{1/2}$, satisfies the three equations in \eqref{F-prime-1}. 

In view of \eqref{eq:kshs}, Eq. \eqref{weingarten-part1} is equivalent to
\begin{equation}
P(x) K_{ext} + Q(x) H = R(x), \label{weingarten-part3}
\end{equation}
where $x=(x_1,x_2,x_3=s)$ and 
\begin{equation}
\begin{array}{rcl}
P(x) &=& (e^{2B}-1)\ce^2(s)-B'\se(2s)-\varphi(x)\se^2(s), \vspace{1ex}\\
Q(x) &=& \epsilon(e^{2B}-1)\se(2s)+2B'\ce(2s)+\varphi(x)\se(2s), \vspace{1ex}\\
R(x) &=& \varphi(x)\ce^2(s)-\epsilon^2(e^{2B}-1)\se^2(s)-\epsilon B' \se(2s).
\end{array}
\end{equation}
Since
$$
\begin{array}{rcl}
4e^{-2B(s)} &=& \overline{P} \ce(2s) + \overline{Q} \se(2s) + \Lambda \se^2(s) + \overline{R},\vspace{1ex} \\
4B'(s) e^{-2B(s)} &=& \epsilon \overline{P} \se(2s) - \overline{Q} \ce(2s)-\Lambda \se(s) \ce(s),
\end{array}
$$
with  $\overline{P}=P+\epsilon R$, $\overline{Q}=Q$, $\Lambda = 2(\epsilon^2-1)R$ and $\overline{R}=P-\epsilon R+4$, we have
$$
\begin{array}{rcl}
4e^{-2B(s)}P(x)&=&-2P- \se^2(s) \left(\epsilon(4e^{-2B(s)}+4 - 2\overline{R})-\Lambda + 4 e^{-2B(s)} \varphi(x) \right),\vspace{1ex}\\
4e^{-2B(s)}Q(x)&=&-2Q+ \se(2s) \left(\epsilon(2e^{-4B(s)}+4 - 2\overline{R}) -\Lambda + 4 e^{-2B(s)} \varphi(x) \right),\vspace{1ex}\\
4e^{-2B(s)}R(x)&=&-2R+ \ce^2(s) \left(\epsilon(2e^{-4B(s)}+4 - 2\overline{R})-\Lambda + 4 e^{-2B(s)} \varphi(x) \right).
\end{array}
$$
The above equations, together with Eqs. \eqref{linear-weingarten-F} and \eqref{weingarten-part3}, imply that
$$
\left(\epsilon(4e^{-2B(s)}+4 - 2 \overline{R}) - \Lambda + 4 e^{-2B(s)} \varphi(x) \right) \left(\se^2(s)K_{ext} - \se(2s) H + \ce^2(s) \right)=0.
$$
We claim that the function 
$$\epsilon(4e^{-2B(s)}+4 - 2\overline{R}) - \Lambda + 4 e^{-2B(s)} \varphi(x)$$
 is identically zero. Let us suppose, by contradiction, that there is a point $x_0$ where such function is nonzero. Then $\se^2(s)K_{ext} - \se(2s) H + \ce^2(s)$ is identically zero in an open neighbourhood $\Omega$ of $x_0$,  hence in $\Omega$ we have
$$
\left\{
\begin{array}{rcl}
\se^2(s)K_{ext} - \se(2s) H &=& - \ce^2(s), \\
P K_{ext} + Q H &=& R.
\end{array}
\right.
$$ 
Since $K_{ext}$ and $H$ depend only on $(x_1, x_2)$, the determinant $Q \se^2(s) + P \se(2s)$ must be identically zero in $\Omega$, which is a contradiction. 
Therefore 
$$\varphi(x) = \dfrac{\epsilon}{2} \left( (\overline{R}-2)e^{2B(s)} - 2 \right) + \dfrac{\Lambda e^{2B(s)}}{4}$$
 for all $x\in M^3$. We conclude that $\varphi_i=0$, $1\leq i\leq 2$, and $\varphi_3= \epsilon (\overline{R}-2)e^{2B(s)}B'(s)$. On the other hand, the function $\theta$ introduced in \eqref{eq:theta} satisfies
$$
\begin{array}{rcl}
2e^{-2B}\theta &=& 2\epsilon (1 - e^{-2B}) + 2e^{-2B} \varphi +  2 e^{-2B} (B'' - 2(B')^2), \\
&=& 2\epsilon (1 - e^{-2B}) - \epsilon(2e^{-2B}+2 - \overline{R}) + \dfrac{\Lambda}{2} - (e^{-2B})'', \\
&=& 2\epsilon (1 - e^{-2B}) - \epsilon(2e^{-2B}+2 - \overline{R})+\dfrac{\Lambda}{2}+\left( 4 \epsilon e^{-2B} - \epsilon \overline{R} - \dfrac{\Lambda}{2} \right), \\
&=& 0.
\end{array}  
$$
Therefore, $\theta=0$. Since
$$
\begin{array}{rcl}
2B'(\epsilon+\varphi) &=& B' \left( 2 \epsilon + \epsilon  \left( (\overline{R}-2)e^{2B} - 2 \right) + \dfrac{\Lambda e^{2B}}{2} \right), \\
&=& \epsilon (\overline{R}-2)e^{2B}B' + \dfrac{\Lambda e^{2B}B'}{2},
\end{array}
$$
we conclude that the third equation in \eqref{F-prime-1} is satisfied.
\qed

\section{Cyclic conformally flat hypersurfaces}

Let $F\colon M^3\to \R^4$ be a conformally flat hypersurface with three distinct principal curvatures $\lambda_1, \lambda_2, \lambda_3$ and corresponding unit principal vector fields $e_1, e_2$ and $e_3$, respectively. E. Cartan proved (see \cite{la}, p. 84) that the conformal flatness of $M^3$ is
equivalent to the relations
\begin{equation} \label{eq:uno}
\left<\nabla_{e_i}e_j,e_k\right>=0
\end{equation}
and
\begin{equation} \label{eq:dos}
(\lambda_j-\lambda_k)e_i(\lambda_i)+(\lambda_i-\lambda_k)e_i(\lambda_j)+ (\lambda_j-\lambda_i)e_i(\lambda_k)=0,
\end{equation}
for all  $1\leq i\neq j\neq k\neq i\leq 3$. It follows from Codazzi's equation and (\ref{eq:uno}) that
\begin{equation}\label{eq:tres}
\nabla_{e_i}e_i=\sum_{j\neq i}(\lambda_i-\lambda_j)^{-1}e_j(\lambda_i)e_j.
\end{equation}
\begin{proposition}\label{prop:equiv}
The following assertions are equivalent:
\begin{itemize}
\item [(i)] The integral curves of $e_1$ are extrinsic circles;
\item [(ii)] The functions 
$\rho_j=\frac{e_j(\lambda_1)}{\lambda_1-\lambda_j}$, $2\leq j\leq 3$, 
satisfy $e_1(\rho_j)=0$. 
\item [(iii)] The relation 
$$
(\lambda_1-\lambda_j)e_je_1(\lambda_1)=2e_1(\lambda_1)e_j(\lambda_1)
$$
holds for $2\leq j\leq 3$.
\item [(iv)] The image by $F$ of each  integral curve $\sigma$ of $e_1$ is contained in a two-dimensional sphere whose normal spaces in $\R^4$ along $F(\sigma)$ are spanned by 
(the restrictions to $F(\sigma)$ of) the vector fields $F_*e_2$ and $F_*e_3$. 
\item [(v)]  The image by $F$ of each leaf of the distribution spanned by $e_2$ and $e_3$ 
is contained in a hypersphere (or affine hyperplane) of $\R^4$;
\end{itemize}
\end{proposition}
\proof The integral curves of $e_1$ are extrinsic circles if and only if 
$$
\left<\nabla_{e_1}\nabla_{e_1}e_1, e_j\right>=0,\;\;\;2\leq j\leq 3.
$$
Using \eqref{eq:uno} and \eqref{eq:tres} we obtain 
\begin{eqnarray*}\left<\nabla_{e_1}\nabla_{e_1}e_1, e_j\right>
&=&e_1\left<\nabla_{e_1}e_1, e_j\right>-\left<\nabla_{e_1}e_1, \nabla_{e_1}e_j\right>\\
&=&e_1(\rho_j),\;\;\;2\leq j\leq 3,
\end{eqnarray*}
hence $(i)$ and $(ii)$ are equivalent.
The equation $e_1(\rho_j)=0$ can be written as 
$$
e_1(\lambda_1-\lambda_j)e_j(\lambda_1)=(\lambda_1-\lambda_j)e_1e_j(\lambda_1), \;\;\;2\leq j\leq 3.
$$
We have
\begin{eqnarray*}
e_je_1(\lambda_1)&=&e_1e_j(\lambda_1)+[e_j, e_1](\lambda_1)\\
&=& e_1e_j(\lambda_1)+(\nabla_{e_j}e_1)(\lambda_1)-(\nabla_{e_1}e_j)(\lambda_1)\\
&=& e_1e_j(\lambda_1)+\left<\nabla_{e_j}e_1,e_j\right>e_j(\lambda_1)-\left<\nabla_{e_1}e_j,e_1\right>e_1(\lambda_1)\\
&=& e_1e_j(\lambda_1)-\frac{e_1(\lambda_j)}{\lambda_j-\lambda_1}e_j(\lambda_1)+\frac{e_j(\lambda_1)}{\lambda_1-\lambda_j}e_1(\lambda_1).
\end{eqnarray*}
Thus, for $2\leq j\leq 3$,  the equation $e_1(\rho_j)=0$ reduces to the relation in item $(iii)$.

Now, for $2\leq j\leq 3$, using (\ref{eq:tres}) we obtain
$$\tilde \nabla_{e_1}F_*e_j= F_*\nabla_{e_1}e_j=-\left<\nabla_{e_1}e_1,e_j\right>F_*e_1=-\rho_jF_*e_1=-\left<F_*e_j, \xi\right>F_*e_1,$$
where  $\xi=\rho_2F_*e_2+\rho_3F_*e_3.$
The equivalence between the assertions in items $(ii)$ and $(iv)$ follows.

Finally, we prove the equivalence between the assertions in items $(iii)$ and $(v)$.
First notice that the normal spaces of the restriction $f_\sigma$ of $f$ to a leaf 
$\sigma$ of the distribution spanned by $e_2$ and $e_3$ are spanned by the restrictions 
to $f(\sigma)$ of $f_*e_1$ and the unit normal vector field $N$ to $F$. 
Since $e_1, e_2$ and $e_3$ are principal directions of $F$, and in view of (\ref{eq:uno}), it follows that $F_\sigma$ has flat normal bundle, with the restrictions of $e_2$ and $e_3$ to $\sigma$ as an orthonormal diagonalizing tangent frame and corresponding principal normal vector fields
\begin{eqnarray*}\eta_j&=&\left<\nabla_{e_j}e_j, e_1\right>F_*e_1+\lambda_jN\\
&=&\frac{e_1(\lambda_j)}{\lambda_j-\lambda_1}F_*e_1+\lambda_jN, \;\;\;2\leq j\leq 3.
\end{eqnarray*}
Using (\ref{eq:dos}) we obtain
\begin{eqnarray*}\eta_2-\eta_3&=&\left(\frac{e_1(\lambda_2)}{\lambda_2-\lambda_1}-\frac{e_1(\lambda_3)}{\lambda_3-\lambda_1}\right)F_*e_1+(\lambda_2-\lambda_3)N\\
&=&(\lambda_2-\lambda_3)(\mu F_*e_1+N),
\end{eqnarray*}
where
$$\mu=\frac{e_1(\lambda_1)}{(\lambda_2-\lambda_1)(\lambda_3-\lambda_1)}.$$
Thus 
\begin{equation}\label{eq:zeta}
\zeta=F_*e_1-\mu N
\end{equation}
is an umbilical normal vector field to $F|_\sigma$, for it is orthogonal to $\eta_2-\eta_3$, and the assertion in item $(v)$ is equivalent to $\zeta$ being parallel with respect to the normal connection of $F|_\sigma$. The latter is, in turn, equivalent to $e_2(\mu)=0=e_3(\mu)$.

Notice that $e_j(\mu)=0$, for $2\leq j\leq 3$, is equivalent to
$$
(\lambda_j-\lambda_1)(\lambda_k-\lambda_1)e_je_1(\lambda_1)
=e_1(\lambda_1)(e_j(\lambda_j-\lambda_1)(\lambda_k-\lambda_1)+e_j(\lambda_k-\lambda_1)(\lambda_j-\lambda_1)),
$$
for $2\leq k\neq j\leq 3$. Using (\ref{eq:dos}), the expression between brackets on the right-hand-side is equal to
$$(\lambda_j-\lambda_k)e_j(\lambda_1)-e_j(\lambda_1)(\lambda_k+\lambda_j-2\lambda_1)=2e_j(\lambda_1)(\lambda_1-\lambda_k).$$
Hence, the equation $e_j(\mu)=0$,  $2\leq j\leq 3$, reduces to the relation in item $(iii)$. \vspace{1ex}\qed

\section{Proof of Theorem \ref{thm:main} and Corollary \ref{cor:ckilling}}

To prove Theorem \ref{thm:main}, let $f\colon M^{2}\times I\to \mathbb{Q}_\epsilon^3\times \R\subset \R_\mu^{4+|\epsilon|}$ be given by (\ref{eq:constantangle}) in terms of an umbilic-free linear Weingarten surface $h\colon\, M^{2}\to \mathbb{Q}_\epsilon^3$ and the function $a\colon I\to \mathbb{R}$ given by \eqref{function-a}. By Theorem \ref{thm:main2},  the metric induced by $f$ on the subset $M^3\subset M^{2}\times I$ of its regular points is conformally flat and  $f$ is a  hypersurface with three distinct principal curvatures such that the tangent component of  $\partial/\partial t$ is a principal direction of $f$ at any point. 

 To complete the proof of the direct statement, it suffices to argue that $\Phi \circ f\colon M^3 \to \mathbb{R}^4$ is a cyclic conformally flat hypersurface, for the composition with an inversion in $\mathbb{R}^4$ clearly preserves both properties. 
 
 We must thus prove that, for each $x\in M^2$, the curve $\gamma\colon I\to M^3$ given by $\gamma(s)=(x,s)$ is a curvature line of $\Phi \circ f$, as well as an extrinsic circle, or a geodesic, of $M^3$. First notice that if $\bar\gamma \colon I\to \mathbb{Q}_\epsilon^3\times \R$ is given by $\bar\gamma= f\circ \gamma$, then $\bar\gamma(I)$ is contained in the vertical cylinder $\beta(\R)\times \R$ in $\mathbb{Q}_\epsilon^3\times \R$ over the geodesic $\beta$ of $\mathbb{Q}_\epsilon^3$ normal to $g$ at $x$, which intersects $f(M)$ orthogonally along $\bar\gamma(I)$. 
 
 We argue separately for the cases $\epsilon=1$, $\epsilon=-1$ and $\epsilon=0$. If $\epsilon =1$, then the image of the vertical cylinder  $\beta(\R)\times \R$ under $\Phi$ is a two-dimensional subspace of $\mathbb{R}^4$ that intersects $\Phi(f(M))$ orthogonally along $\Phi(\bar\gamma(I))$, for $\Phi$ is conformal. 
Thus $\gamma$ is a curvature line of $\Phi \circ f$ and also a geodesic of $M^3$ (see Proposition~$9$ of \cite{to1}).

  If $\epsilon =-1$, then the image of the vertical cylinder $\beta(\R)\times \R$ under
    $\Phi\colon \mathbb{H}^3\times \R\to \R^{4}\setminus \R^2$ is a two-dimensional sphere centered at the subspace $\mathbb{R}^2\subset \mathbb{R}^4$, which intersects $\Phi(f(M))$ orthogonally along $\Phi(\bar\gamma(I))$. Therefore, in this case the curve $\gamma$ is a curvature line of $\Phi \circ f$ that is an extrinsic circle of $M^3$ (see again Proposition~$9$ of \cite{to1}). The case $\epsilon=0$ is similar and easier.
   
    To prove the converse statement, let $F\colon M^3 \to \R^4$  be a cyclic conformally flat hypersurface. By Proposition \ref{prop:equiv}, since the integral curves of $e_1$ are extrinsic circles, one has a family $\mathcal{F}$ of hyperspheres 
(or affine hyperplanes) that contain the images by $F$ of the leaves of the distribution spanned by $e_2$ and $e_3$, and a family $\mathcal{G}$ of two-dimensional spheres 
(or affine subspaces) that contain the images 
by $F$ of the integral curves of $e_1$, with the property that each element of the former is 
orthogonal to every element of the latter, and conversely. 
By Lemma $6$ of \cite{to1}, there exists an inversion $\mathcal{I}$ in $\R^4$ that takes the families $\mathcal{F}$ and $\mathcal{G}$, respectively, into families of hyperspheres 
(or affine hyperplanes) and two-dimensional spheres (or affine subspaces) of one of the following types:
\begin{itemize}
\item[(i)] a family of parallel affine hyperplanes and a family of orthogonal affine subspaces;
\item[(ii)] a family of concentric hyperspheres and a family of affine subspaces through their common center;
\item[(iii)]  a family of affine hyperplanes intersecting along a two-dimensional affine subspace 
and a family of two-dimensional spheres centered at that affine subspace;
\item[(iv)] a family of hyperspheres whose centers lie in a straight line and a family of two-dimensional 
affine subspaces intersecting along that straight line;
\item[(v)] a family of affine hyperplanes intersecting along a straight line and a family of two-dimensional 
spheres centered at that straight line.
\end{itemize}

 Let $\zeta$ be the vector field given by \eqref{eq:zeta}. In case $(i)$, the vector field $\mathcal{I}_*\zeta$ is collinear with the constant unit vector field $e_4$ normal to the family $\mathcal{I}(\mathcal{F})$ of affine hyperplanes, thus the tangent component of $e_4$ is collinear with $\mathcal{I}_*F_*e_1$, and hence is a principal direction of $\tilde f=\mathcal{I}\circ F$. In terms of the orthogonal decomposition  $\mathbb{R}^4=\mathbb{R}^3 \times \mathbb{R}$, with $e_4$ spanning the factor $\mathbb{R}$, we can write $\tilde f=\Phi\circ f$, and hence $F=\mathcal{I}\circ \Phi\circ f$, where $\Phi\colon \mathbb{R}^3 \times \mathbb{R}\to  \mathbb{R}^4$ is the standard isometry and $f\colon M^3\to \mathbb{R}^3 \times \mathbb{R}$ is a conformally flat hypersurface with three distinct principal curvatures, having the property that the tangent component of the vector field $\partial/\partial t=e_4$ is a principal direction at any point. 

In case $(ii)$, the vector field $\mathcal{I}_*\zeta$ is collinear with the radial vector field $\mathcal{R}$ along $\tilde f=\mathcal{I}\circ F$,  thus the tangent component of $\mathcal{R}$ along $\tilde f$ is collinear with $\mathcal{I}_*f_*e_1$, and hence is a principal direction of $\tilde f$.
In other words, $F=\mathcal{I}\circ \tilde f$, where $\tilde f\colon M^3\to \R^4$ is a conformally flat hypersurface with three distinct principal curvatures having the property that the tangent component of the radial vector field $\mathcal{R}$ along $\tilde f$, that is, of the position vector field of $\tilde f$,  is a principal direction of $\tilde f$. 
It follows that $\tilde f=\Phi \circ  f$, and hence $F=\mathcal{I}\circ \Phi\circ f$, where $\Phi\colon \mathbb{S}^3\times \R\to \R^4\setminus \{0\}$ is the conformal diffeomorphism given by
$\Psi(x,t)=e^tx$
and $f\colon M^3\to \mathbb{S}^3\times \R$ is a conformally flat hypersurface with three distinct principal curvatures having the property that the tangent component of the unit vector field $\frac{\partial}{\partial t}$ is a principal direction at any point, for the vector fields $\frac{\partial}{\partial t}$ and $\mathcal{R}$ are $\Phi$-related, that is,
$\Phi_*(x,t)\frac{\partial}{\partial t}=\mathcal{R}(\Phi(x,t))$.

  In case $(iii)$ we assume that the affine subspace in the intersection of all affine hyperplanes of the family  $\mathcal{I}(\mathcal{F})$ is, say, 
  the subspace
$\{(y_0, y_1, y_2, y_{3})\,:\, y_0=0=y_1\}$. 
Then the vector field $\mathcal{I}_*\zeta$ is collinear along $\tilde f=\mathcal{I}\circ F$ with the Killing vector field $\mathcal{K}$ in $\R^4$ given by
$ \mathcal{K}(y_0, y_1, y_2, y_{3})=(0, -y_{3}, y_2).$
Thus the tangent component of $\mathcal{K}$ along $\tilde f$ is collinear with $\tilde{f}_*e_1=\mathcal{I}_*F_*e_1$, and hence is a principal direction of $\tilde f$. Now notice that  $\mathcal{K}$ and the unit vector field $\partial/\partial t$ tangent to the factor $\mathbb{R}$ in $\mathbb{H}^3\times \R$ are $\Phi$-related, where $\Phi\colon \mathbb{H}^3\times \R\to \R^4\setminus \R^2$ be the  conformal covering map given by \eqref{tpsi}, that is, $\Phi_*(x,t)\frac{\partial}{\partial t}=\mathcal{K}(\Phi(x,t))$. It follows that $\tilde f=\Phi \circ f$, and hence $F=\mathcal{I}\circ \Phi\circ f$,  where $f\colon M^3\to \mathbb{H}^3\times \R$ is a conformally flat hypersurface with  three distinct principal curvatures having the property that the tangent component of the unit vector field $\frac{\partial}{\partial t}$ is a principal direction at any point.

 In all three cases above, it follows from Theorem \ref{thm:main2} that $f\colon M^3\to \mathbb{Q}_\epsilon \times \mathbb{R}$ is given by (\ref{eq:constantangle}) in terms of a linear Weingarten surface $h\colon\, M^{2}\to \mathbb{Q}_\epsilon^3$ and a smooth function $a\colon I\to \mathbb{R}$ given by \eqref{function-a}. 

We now argue that cases $(iv)$ and $(v)$ can not occur. 
Let $(y_1, \ldots, y_4)$ be standard coordinates on $\R^4$ and let  $\Psi\colon  \R^{4}\setminus \R \to 
\mathbb{H}^2 \times \mathbb{S}^{2} \subset \R^{3}_1 \times \R^{3}$ be the conformal diffeomorphism, with $\R=\{(y_1, \ldots, y_4) \in \R^4\,:\, y_2=y_3=y_4=0\}$, given by
$$ \Psi(y_1,y_2,y_3, y_4) = 
\frac{1}{\sqrt{y_2^2+y_3^2+y_4^2}} \left( e_0 + y_1e_1
+ \left( \sum_{i=1}^{4} y_i^2 \right) e_2, (y_2, y_3, y_4) \right), $$
where $e_0, e_1, e_2$ is a 
pseudo-orthonormal basis of $\R_1^3$ with $\left<e_0, e_0\right>=0=\left<e_2, e_2\right>$, $\left<e_0, e_2\right>=-1/2$ and $\left<e_1,e_j\right>=\delta_{1j}$, $0\leq j\leq 2$. 

If either $(iv)$ or $(v)$ holds, then
$ f=\Psi\circ \mathcal{I}\circ F\colon M^3\to \mathbb{H}^2 \times \mathbb{S}^{2}$ maps each integral curve of $e_1$ into a slice $\mathbb{H}^2\times \{x\}$ or $\{x\}\times \mathbb{S}^2$ of $\mathbb{H}^2 \times \mathbb{S}^{2}$, respectively.  
 In the former case, $f(x,s)=(a(s), h(x,s))$ for some smooth maps $a\colon I\to \R$ and $h\colon M^3\to \mathbb{S}^2$. Since $f_*\partial_s$ is orthogonal to
 $f_*X$ for any $X\in T_xM^2$, it follows that $h$ does not depend on $s$. But then $\Phi \circ f=\mathcal{I}\circ F$ would be a rotation hypersurface over the plane curve $s\mapsto \Phi(a(s), h(x))$, for a fixed $x\in M^2$. Therefore $\mathcal{I}\circ F$, and hence also $F$, would have only two distinct principal curvatures, a contradiction. Arguing in a similar way also rules out case $(v)$. \vspace{1ex} \qed
 
 \noindent \emph{Proof of  Corollary \ref{cor:ckilling}:} If $F\colon M^3\to \R^4$ is a cyclic conformally flat hypersurface, by Theorem~\ref{thm:main} it is given by $F=\mathcal{I} \circ \Phi \circ f\colon M^3 \to \mathbb{R}^4$ in terms of a hypersurface 
 $f\colon M^3 \to \mathbb{Q}_\epsilon^3\times \mathbb{R}$  as in Theorem~\ref{thm:main2}, where $\mathcal{I}$ is either the identity map or an inversion with respect to a hypersphere  in $\R^4$ and $\Phi$ denotes either the conformal diffeomorphism 
 $\Phi\colon \mathbb{Q}_\epsilon^3\times \R\to \R^{4}\setminus \{0\}$ if $\epsilon=1$, the conformal covering map $\Phi\colon \mathbb{Q}_{\epsilon}^3\times \R\to \R^{4}\setminus \R^2$ if $\epsilon=-1$ or the isometry $\Phi\colon \mathbb{Q}_\epsilon^3\times \R\to  \R^4$ if $\epsilon=0$. 
 
 Since $\Phi$ is a conformal diffeomorphism, the tangent component of  $\partial/\partial t$ is a principal direction of $f$ at any point and $\partial/\partial t$ is $\Phi$-related to either a constant vector field $\partial_{x_i}$, the radial vector field $\mathcal{R}$ or one of the Killing vector fields $\mathcal{K}_{ij}$ in $\R^4$, according to whether $\epsilon=0$, $\epsilon=1$ or $\epsilon=-1$, respectively, then the tangent component of one of those vector fields is a principal direction of $\tilde f=\Phi\circ f$ at any point. Finally, if $\mathcal{I}$ is an inversion with respect to a hypersphere  in 
$\R^4$, then (a multiple of) the vector field $\partial_{x_i}$ is $\mathcal{I}$-related to $\mathcal{C}_i$, whereas $\mathcal{R}$ is  $\mathcal{I}$-related to (a multiple of) itself. 
 Therefore,  the tangent component of either $\mathcal{C}_i$ or $\mathcal{R}$ is a principal direction of $F=\mathcal{I} \circ \Phi \circ f$ at any point.
 
 Conversely, assume that  $F\colon M^3\to \R^4$ is a  conformally flat hypersurface with three distinct principal curvatures such that the tangent component of one of the conformal Killing vector fields $\partial_{x_i}$,  $\mathcal{R}$, $\mathcal{K}_{ij}$ or 
 $\mathcal{C}_i$ is a principal direction of $\tilde f=\Phi\circ f$ at any point. We argue for  $\mathcal{C}_i$, the other cases being similar. Since (a multiple of) the vector field $\partial_{x_i}$ is $\mathcal{I}$-related to $\mathcal{C}_i$, it follows that the tangent component of  $\partial_{x_i}$ is a principal direction of $\tilde f=\mathcal{I}\circ F$ at any point. Let $\Phi\colon \R^3\times\R\to  \R^4$ be the isometry given by the orthogonal decomposition of $\R^4$ determined by $\partial_{x_i}$. Then $F=\mathcal{I}\circ \Phi\circ f$, where $f\colon M^3\to \R^3\times\R$ has the property that the tangent component of the unit vector field  $\partial/\partial t$ tangent to $\R$ is a principal direction of $f$ at any point. Thus $F$ is a cyclic conformally flat hypersurface by Theorem \ref{thm:main}. \qed

{\renewcommand{\baselinestretch}{1}
\hspace*{-30ex}\begin{tabbing}
\indent \= Universidade de Brasilia  \hspace{11.5ex} Universidade de S\~ao Paulo \\
\>  Campus Univ. Darcy Ribeiro \hspace{6.5ex}
Av. Trabalhador S\~ao-Carlense 400 \\
\> 70910-900 --- Brasilia -- DF
\hspace{8.5ex} 13560-970 --- S\~ao Carlos -- SP \\
\> Brazil\hspace{31ex} Brazil\\
\> joaopsantos@unb.br  \hspace{16ex}
tojeiro@icmc.usp.br
\end{tabbing}}

\begin{thebibliography}{lbllll}

\bibitem[{\bf CT$_1$}]{ct1}  Canevari, S. and Tojeiro, R.,  \emph{Hypersurfaces of two space forms and conformally flat hypersurfaces}. 
Ann. Mat. Pura  Appl. 197 (2018), 1--20.

\bibitem[{\bf CT$_2$}]{ct2}  Canevari, S. and Tojeiro, R.,  \emph{The Ribaucour transformation for  hypersurfaces of two space forms and conformally flat hypersurfaces}. 
Bull. Braz. Math. Soc. 49 (2018), 593--613. 


\bibitem[{\bf Ca}]{ca} Cartan, E., {\it La d\'eformation des hypersurfaces dans l'espace conforme r\'eel a $n\geq 5$ dimensions.}
Bull. Soc. Math. France {\bf 45} (1917), 57--121.

\bibitem[{\bf DT$_1$}]{dt1} Dajczer, M. and Tojeiro, R., \emph{ On a class of submanifolds carrying an extrinsic umbilic foliation}. Israel J. Math. 125 (2001), 203--220.

\bibitem[{\bf DT$_2$}]{dt2} Dajczer, M. and Tojeiro, R., ``Submanifold Theory beyond an introduction", Springer, New York, 2019, Universitext Series.

\bibitem[{\bf H-J}]{hj} Hertrich--Jeromin, U., {\it On Conformally Flat Hypersurfaces and Guichard's Nets.\/}  Beitr. Alg. Geom. (1994), 315--331.

\bibitem[{\bf H-J}$_2$]{hj2} Hertrich--Jeromin, U.,
``Introduction to M\"obius differential geometry". 
London Mathematical Society Lecture Note Series, 300. 
Cambridge University Press, Cambridge, 2003.

\bibitem[{\bf H-JS}]{hjs} Hertrich--Jeromin, U. and Suyama, Y., \emph{Conformally Flat Hypersurfaces with Cyclic  Guichard Net}, 
Int. J. Math. 18 (2007), 301--329. 

\bibitem[{\bf H-JS$_2$}]{hjs2} Hertrich--Jeromin, U. and Suyama, Y.,  
\emph{Conformally flat hypersurfaces with Bianchi-type  Guichard net},
Osaka J. Math. 50   (2013), 1--30.

\bibitem[{\bf H-JSUY}]{hjsuy} [Hertrich--Jeromin, U., Suyama, Y., Umehara, M. and Yamada, K., \emph{A duality for conformally flat hypersurfaces}, Beitr. Alg. Geom. 56 (2018), 655--676.

\bibitem[{\bf La}]{la} Lafontaine, J., {\it Conformal geometry from the Riemannian viewpoint.\/} Aspects of Mathematics, E 12, Vieweg, Braunschweig,
1988.

\bibitem[{\bf ST}]{st} dos Santos, J. P. and Tenenblat, K., \emph{The symmetry group of 
Lam\'e system and the associated Guichard nets for conformally flat hypersurfaces},
SIGMA 9 (2013), p. 33.

\bibitem[{\bf Su1}]{su1}  Suyama, Y., \emph{Conformally flat hypersurfaces in Euclidean $4$-space $II$},
Osaka J. Math. 42 (2005) 573--598.


\bibitem[{\bf Su2}]{su2}  Suyama, Y., \emph{A classification and a non-existence theorem for conformally flat hypersurfaces in Euclidean 4-space}, Int. J. Math. 16 (2005) 53--85.


\bibitem[{\bf To$_1$}]{to1} Tojeiro, R, \emph{Conformal immersions of warped products}, Geom. Dedicata 128  (2007), 17--31.

\bibitem[{\bf To$_2$}]{to2} Tojeiro, R, \emph{On a class of hypersurfaces in $\mathbb{S}^n\times \R$ and $\mathbb{H}^n\times \R$},  Bull. Braz. Math. Soc. 41  (2010), 199--209.


\end{thebibliography}
\end{document}